\documentclass{article}
\usepackage{arxiv}

\usepackage{amsmath,amssymb}
\usepackage{changepage}
\usepackage{textcomp,marvosym}
\usepackage[nopatch=eqnum]{microtype}
\usepackage{nameref,hyperref}
\usepackage[right]{lineno}
\DisableLigatures[f]{encoding = *, family = * }

\usepackage{makecell}
\usepackage{graphicx}
\usepackage{subcaption}

\usepackage[table]{xcolor}

\usepackage{array}

\newcolumntype{+}{!{\vrule width 2pt}}

\newlength\savedwidth

\usepackage{tabularx}
\usepackage{array}

\usepackage{lastpage,fancyhdr,graphicx}
\usepackage{epstopdf}
\usepackage{ifthen}

\usepackage{booktabs}
\usepackage{comment}
\usepackage{mathtools,amsfonts,cases}
\usepackage{lipsum}

\usepackage{tikz}
\usepackage{import}
\usepackage{xcolor}
\usepackage{comment}
\usepackage{cancel}
\pgfkeys{/pgf/number format/.cd,1000 sep={}}

\usepackage{rotating}
\usepackage{multirow}

\usepackage{setspace}

\usepackage{graphicx}
\graphicspath{{figures/}}

\usepackage[numbers]{natbib}

\title{Second-Order IMEX Time-Stepping Methods for Efficient Bidomain Multi-Electrode Array Simulations of Cardiac Stem Cell Monolayers}

\author{\href{https://orcid.org/0009-0009-8203-2164}{\includegraphics[scale=0.06]{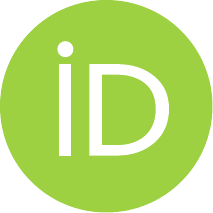}\hspace{1mm}\textbf{S.~ Tonali}}\thanks{\texttt{sofia.tonali01@universitadipavia.it}
Department of Mathematics, University of Pavia, Italy;
Euler Institute, Università della Svizzera Italiana, Lugano, Switzerland
}
\\
\href{https://orcid.org/0000-0002-7271-6943}{\includegraphics[scale=0.06]{orcid.pdf}\hspace{1mm}\textbf{S.~Botti}}\thanks{\texttt{{sofia.botti@polimi.it}}, MOX, Department of Mathematics, Politecnico di Milano, Italy;
Euler Institute, Università della Svizzera Italiana, Lugano, Switzerland}
\\
\href{https://orcid.org/0000-0002-3014-4668}{\includegraphics[scale=0.06]{orcid.pdf}\hspace{1mm}\textbf{L.~F.~Pavarino}}\thanks{\texttt{{luca.pavarino@unipv.it}}, Department of Mathematics, University of Pavia, Italy}
}

\renewcommand{\shorttitle}{\textit{arXiv} Template}

\newcommand{\MRMS}{\ensuremath{\mathrm{MRMS}}}
\newcommand{\Nel}{\ensuremath{\mathrm{N_{el}\,}}}
\newcommand{\ooc}[1]{\ensuremath{p_{\mathrm{#1}}}}

\usepackage{xparse}
\NewDocumentCommand{\err}{o o}{%
  \ensuremath{%
    \mathcal{E}%
    \IfValueT{#1}{_{\mathnormal{#1}}}%
    \IfValueT{#2}{^{\mathnormal{#2}}}%
  }%
}

\NewDocumentCommand{\FP}{o}{%
  \ensuremath{%
    U_{\mathrm{FP}}%
    \IfValueT{#1}{^{#1}}%
  }%
}

\usepackage{mathrsfs}

\begin{document}
\maketitle

\onehalfspacing
\bibliographystyle{plainnat}

\begin{abstract}
Multi-electrode Arrays (MEAs) enable tissue-level electrophysiological studies of human induced pluripotent stem cell-derived cardiomyocytes (hiPSC-CMs) by recording extracellular field potentials. Recent computational models have improved MEA simulations by coupling detailed electrode descriptions with the Bidomain framework, a parabolic-elliptic system of nonlinear PDEs coupled with a stiff ionic model. The standard numerical strategy relies on operator splitting techniques that decouple the PDE and ODE components. In most implementations, the ionic subsystem is treated explicitly while the diffusive operator is handled implicitly, resulting in a first-order implicit-explicit (IMEX) time discretization. Although computationally convenient, this approach limits temporal accuracy and may reduce efficiency in large-scale simulations.

In this work, we investigate higher-order IMEX Runge–Kutta schemes within a Strang-based operator splitting, specifically tailored for the MEA model of hiPSC-CMs monolayers. First- and second-order schemes are compared in terms of computational cost and global error against a high-fidelity reference solution obtained with a very small time step. For a fixed computational cost, the second-order schemes achieve errors that are 2–3 orders of magnitude smaller than those of the first-order method. These results demonstrate that higher-order IMEX integration significantly improves the accuracy-to-cost ratio of MEA simulations, providing a practical and reliable approach for large-scale electrophysiological studies.
\end{abstract}

\section{Introduction}

Numerical simulation of cardiac electrophysiology plays a central role in understanding the mechanisms underlying electrical propagation in the heart and in supporting the development of in-vitro experimental platforms. A major breakthrough in this field has been the discovery of human induced pluripotent stem cells (hiPSCs) \cite{Yamanaka2007, Yamanaka2006}. These cells, obtained by reprogramming adult somatic cells, can be differentiated into various lineages, including cardiomyocytes (CMs). Since they express the main cardiac ion channels and exhibit electrophysiological responses comparable to adult human cells, hiPSC-CMs have become invaluable tools for drug safety pharmacology and disease modelling both in vitro (\cite{karakikes_invitro,wu2024clinical}) and in silico (\cite{Tveito_2,Tveito,simone2026novel}).

From a mathematical perspective, existing hiPSC-CM tissue simulations rely on the Bidomain model \cite{tung1978}, a reaction-diffusion system coupled with a stiff system of ordinary differential equations representing ionic currents \cite{Botti2024,Kernik,Koivumaki_hiPSC,Paci2013,Paci2020,Paci2018}. The strong multiscale nature of this problem presents significant challenges for numerical discretization, requiring methods that balance stability, accuracy and computational efficiency. 
 
A widely adopted approach for simulating these models is using Finite Element Method (FEM) in space and an operator splitting (OS) technique, which decouples the fully coupled nonlinear systems of PDEs into simpler sub-problems, in time. In this context, Implicit--Explicit (IMEX) schemes are particularly effective, allowing diffusive terms to be handled implicitly while treating the nonlinear ionic dynamics explicitly. A comprehensive review of IMEX methods in a general framework can be found in \cite{Russo}. Regarding cardiac modelling, early investigations were conducted in \cite{spiteri2008}, analysing the performance of IMEX Runge--Kutta (IMEX--RK) methods and highlighting their favourable stability properties and computational efficiency. More recent works (\cite{cervi2018high,cervi2019fourth,spiteri2016}) have further explored higher-order OS and IMEX schemes for cardiac simulations, demonstrating their potential to improve accuracy without significantly increasing computational costs.

In more recent years, Multi-Electrode Arrays (MEA) have emerged as a powerful tool for investigating the electrophysiological behaviour of hiPSC-CMs, enabling non-invasive recordings of extracellular field potential (FP) at the tissue level (\cite{Pickard1979,SpiraMEA}). These experimental advances call for accurate and efficient in silico models capable of reproducing both the intracellular dynamics and the extracellular FP measured by MEA. Thus, to faithfully reproduce MEA recordings, the standard bidomain formulation must be extended to incorporate a detailed description of the electrode--tissue interface.

Although higher-order IMEX schemes have shown promising results in cardiac electrophysiology, their effectiveness within detailed MEA simulations and hiPSC-CM ionic models remains to be further explored. In this work, we investigate the application of higher-order IMEX Runge--Kutta schemes within a detailed MEA simulation framework for hiPSC-CM monolayers. Specifically, we assess the performance of second-order IMEX--RK methods combined with Strang-based operator splitting \cite{Strang} and compare them against the first-order approach. The goal is to assess whether increasing the temporal order provides tangible efficiency gains in a MEA model that includes a detailed electrode description. A significant challenge lies in the complexity of the ionic model: for ventricular-like hiPSC-CMs, we employ the Paci 2020 model \cite{Paci2020}, which involves 23 ionic currents and 3 ionic concentrations, resulting in a highly stiff and computationally demanding system.

The paper is organized as follows. Section \ref{sec:model} presents the MEA model together with the ionic dynamics and the main modelling assumptions. Section \ref{sec:NumericalMethods} describes the numerical methods, including the construction of the reference solution and the error analysis framework used to assess convergence properties. Numerical results are presented in Section \ref{sec:results}. We first analyse the performance of the proposed schemes for the zero-dimensional ionic model to isolate the temporal integration error before introducing spatial complexity, and subsequently consider the full two-dimensional MEA setting. Finally, Section \ref{sec:conclusion} summarizes the main findings of the study, discusses its limitations, and outlines possible directions for future research.

\section{Mathematical Model}
\label{sec:model}
\begin{figure}[t]
    \centering
    \includegraphics[width=\linewidth]{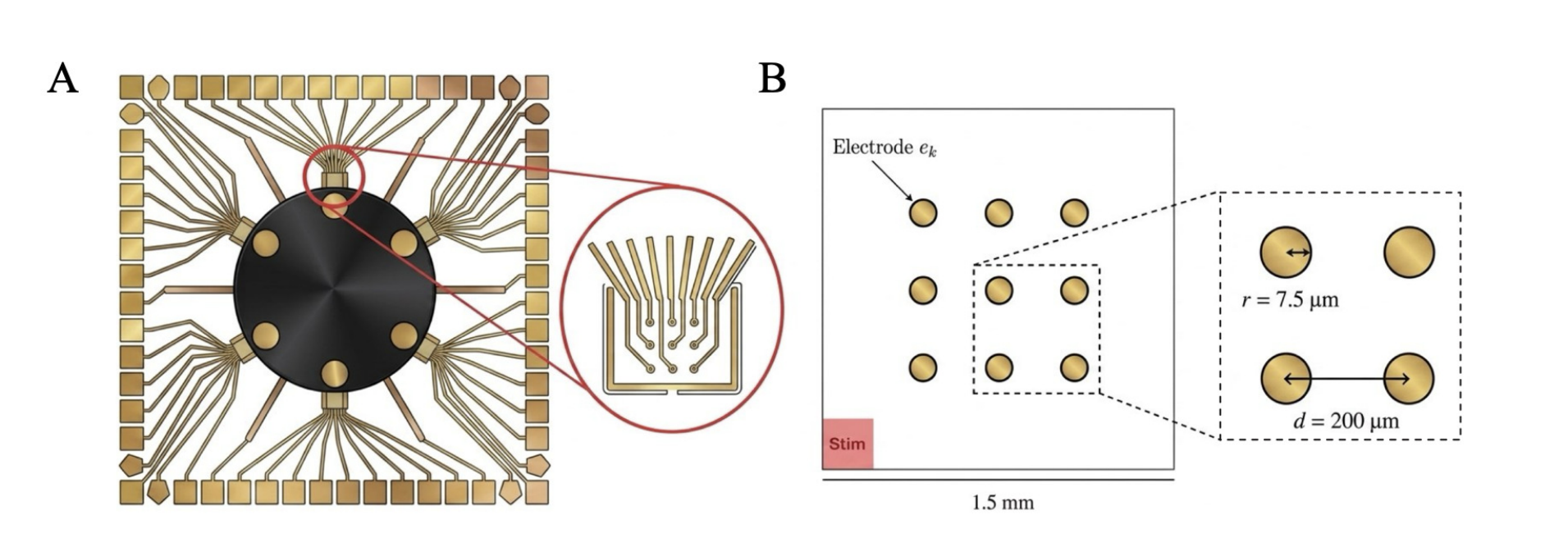}
    \caption{(A) Schematic representation of the experimental MEA setup. Each well contains a 3$\times$3 array of recording electrodes and an internal reference electrode.
    (B) Computational geometry adopted in the MEA model simulations. The domain includes the 3$\times$3 equally spaced electrode configuration. Electrical propagation is initiated by a stimulus given in the bottom-left corner of the domain.}
    \label{MEA_Framework}
\end{figure}

The mathematical model considered in this work is based on the bidomain equations, widely used for the description of the electrical activity of the cardiac tissue. Since in our case MEA devices are considered, to model this setting, the bidomain model is coupled with a model for the recording electrodes as done in \cite{Coudier,Botti_MEA}. 
Our framework reproduces the 60–6 well MEA platform system from Multi Channel Systems (see \cite{MultiChannelSystem} for details), but the method can be applied to other kinds of MEAs. This device consists of six independent wells, each equipped with nine recording electrodes (see Fig. \ref{MEA_Framework}). Six independent experiments can thus be done with identical surrounding conditions at once MEA, but for simplicity, in the present study, simulations are carried out at the single-well level. For simplicity, throughout this work, we will refer to the coupled bidomain–MEA model simply as the MEA model.

\subsection*{The MEA model} 
Let $\Omega \subset \mathbb{R}^2$ denote a monolayer of hiPSC-CMs connected to a MEA device with \Nel electrodes and $(0,T)$ the time interval considered. We use the bidomain equations to describe the electrical activity of the tissue. The bidomain model describes cardiac electrical propagation in a spatial domain as the superposition of two anisotropic conductive media: the intracellular and extracellular compartments (a detailed derivation and more information can be found in \cite{pavarino} and \cite{sundnes}). In a two-dimensional setting, the full MEA model in a parabolic-elliptic (PE) form can be written as:

\begin{subequations}\label{eq:mea_formulation}
\begin{align}
& C_m \frac{\partial v}{\partial t} - \nabla(\mathrm{D_i} \nabla v) - \nabla(\mathrm{D_i} \nabla u_e) + I_{\text{ion}} = I_{\text{app}} && \text{in } \Omega \times (0,T) \label{mea_1} \\[2pt]
& - \nabla(\mathrm{D_i} \nabla v) - \nabla((\mathrm{D_i} + \mathrm{D_e}) \nabla u_e) = \frac{1}{\mathrm{z_{\text{thick}}}} \sum_{e_k} \frac{I_{el}^k}{|e_k|} \chi_{e_k} && \text{in } \Omega \times (0,T) \label{mea_2} \\[2pt]
& \frac{dI_{el}^k}{dt}+\frac{I_{el}^k}{\tau}=\frac{C_{el}}{\tau} \frac{dU^k}{dt}, \quad U^k=\frac{1}{|e_k|} \int_{e_k} u_e \mathrm{d}e_k && \text{in } e_k \times (0,T) \label{mea_3} \\[1pt]
& &&k=1,\dots, \Nel \notag \\[2pt]
& \frac{d\boldsymbol{w}}{dt} - \mathbf{W}(v, \boldsymbol{w}, \boldsymbol{c})= 0, \quad \frac{d\boldsymbol{c}}{dt} - \mathbf{C}(v, \boldsymbol{w}, \boldsymbol{c})= 0 && \text{in } \Omega \times (0,T) \label{ionic_model}
\end{align}
\end{subequations}

 The functions $u_i(\boldsymbol{x},t), u_e(\boldsymbol{x},t)$ are, respectively, the intracellular and the extracellular potential, $v(\boldsymbol{x},t) = u_i(\boldsymbol{x},t) - u_e(\boldsymbol{x},t)$ is the transmembrane potential, and its time course (waveform) over an excitation - recovery cycle is known as action potential (AP). Electrical fluxes in each compartment obey diffusion laws governed by the conductivity tensors $\mathrm{D_i}$ and $\mathrm{D_e}$, for the intra- and extracellular media, respectively. Charge exchange across the membrane depends on the membrane capacitance $C_m$. 
 The functions $\boldsymbol{w}(\boldsymbol{x},t), \boldsymbol{c}(\boldsymbol{x},t)$ are, respectively, the gating variables and ionic concentrations of the ionic model associated with the ionic current $I_{\text{ion}} = I_{\text{ion}}(v, \boldsymbol{w}, \boldsymbol{c})$ in \eqref{mea_1} and the differential systems in \eqref{ionic_model}, described below.
 An externally applied stimulus current, denoted by $I_{\mathrm{app}}$, can be applied on a specific region of the domain for a specific time duration. The surface of the electrode $k$ in contact with the monolayer is denoted with $e_k$, while $\chi_{e_k}$ is its characteristic function, i.e. the function equal to one inside the electrode and zero outside. Finally, $I_{el}^k$ is the recorded current by the $k^{th}$ electrode. To compute $I_{el}^k$ we use the electric model described by the equation \eqref{mea_3}, where $\tau=(R_i+R_{el})C_{el}$, $R_i$ denoting the inner resistance, $R_{el}$ and $C_{el}$ the resistance and
the capacitance of the electrode. In equation \eqref{mea_3} we have also introduced $U^k$ as the average extracellular potential over $e_k$. Since the thickness of the monolayer of cells $\mathrm{z_{thick}}$ is supposed to be very small compared to the other dimensions of the problem, we assume that all the variations along the $z$-direction are negligible compared to the one along $x$- and $y$-direction. Finally, the measured field potential at electrode $k$ is computed as
\begin{equation}
\FP[k] = R_i I_{el}^k.
\end{equation}

\subsection*{Initial and Boundary Conditions} 
The MEA system \eqref{eq:mea_formulation} is completed with the initial conditions:
\begin{equation}
v(\boldsymbol{x},0) = v_0, \quad \boldsymbol{w}(\boldsymbol{x},0) = \boldsymbol{w}_0, \quad \boldsymbol{c}(\boldsymbol{x},0) = \boldsymbol{c}_0,
\end{equation}
and with boundary conditions of mixed type. The intracellular potential, defined as $u_i = v + u_e$, satisfies a no-flux condition on the entire boundary:
\begin{equation}
(\mathrm{D}_i \nabla u_i) \cdot \mathbf{n} = 0 \quad \text{on } \partial\Omega \times (0,T).
\end{equation}
The extracellular potential $u_e$ is grounded on a portion of the boundary $\partial\Omega_D \subseteq \partial\Omega$, while a no-flux condition is imposed on the remaining part $\partial\Omega_N = \partial\Omega \setminus \partial\Omega_D$:
\begin{equation}
\begin{cases}
u_e = 0 & \text{on } \partial\Omega_D \times (0,T), \\
(\mathrm{D}_e \nabla u_e) \cdot \mathbf{n} = 0 & \text{on } \partial\Omega_N \times (0,T),
\end{cases}
\end{equation}
Here $\mathbf{n}$ denotes the outward unit normal vector to $\partial\Omega$.

\subsection*{Ionic model: the Paci 2020 model} 
The ionic current $I_{\text{ion}} = I_{\text{ion}}(v, \boldsymbol{w}, \boldsymbol{c})$ in \eqref{mea_1} and the differential systems in \eqref{ionic_model}, which reproduce the electrophysiological activity of individual cells, including the functions of various transmembrane channels, depends on the considered
ionic model.  While this study focuses on hiPSC-derived tissue, the framework can be easily adapted to adult cardiac tissue by employing suitable ionic models, such as the Luo-Rudy \cite{luo-rudy} or the ten Tusscher \cite{tenTusscher} models. 
In this work, we employ the Paci 2020 model \cite{Paci2020} to reproduce the ventricular-like (VL) phenotype of hiPSC-CMs and capture their fundamental electrophysiological features, including their spontaneous automaticity. The Paci 2020 model accounts for several transmembrane currents and intracellular calcium dynamics, resulting in a system of ODEs that is highly nonlinear and stiff. 
In particular, the vector $\boldsymbol{w}$ collects 23 state variables, each representing a specific gating variable associated with an ionic current. The vector of concentrations $\boldsymbol{c}$ includes three ionic species: intracellular sodium $[\text{Na}]_i$, intracellular calcium $[\text{Ca}]_i$, and the calcium concentration within the sarcoplasmic reticulum $[\text{Ca}]_{\text{SR}}$. Further details regarding the Paci model and its previous formulations can be found in \cite{Paci2013, Paci2020, Paci2018}.

\subsection*{Tissue conductivity}
Electrical propagation in cardiac tissue is strongly influenced by cellular alignment and structural organization, which in mature myocardium lead to anisotropic conduction properties. In general, the intra- and extracellular conductivity tensors $\mathrm{D}_{i,e}(\boldsymbol{x})$ can be expressed in anisotropic form as
\[
\mathrm{D}_{i,e}(\boldsymbol{x}) = \sigma_{i,e}^l \, \boldsymbol{a}_l(\boldsymbol{x}) \boldsymbol{a}_l^\top(\boldsymbol{x})
+ \sigma_{i,e}^t \, \boldsymbol{a}_t(\boldsymbol{x}) \boldsymbol{a}_t^\top(\boldsymbol{x}),
\]
where $\boldsymbol{a}_l$ and $\boldsymbol{a}_t$ denote orthonormal principal directions and $\sigma_{i,e}^l$, $\sigma_{i,e}^t$ the corresponding longitudinal and transverse conductivities.

In hiPSC-CM monolayers, however, cell orientation, gap-junction distribution, and maturation state can vary substantially, making an accurate characterization of directional conductivities difficult (see \cite{Pasqualini},\cite{lodrini2020}). For these reasons we assume an isotropic medium with
\[
\mathrm{D}_{i,e} = \sigma_{i,e}\mathrm{I}.
\]
Following \cite{Tveito_2}, we set $\sigma_i = 0.01$~mS/cm and $\sigma_e = 0.1$~mS/cm.

\subsection{Applied Current}
A stimulation current of magnitude $I_{\mathrm{app}} = 150~\mathrm{\mu A/cm^2}$ was applied to a square region of $3600 ~\mu \text{m}^2$ at the bottom-left corner of the domain for a duration $T_{\mathrm{stim}}=0.51$ ms. To improve the temporal regularity of the stimulus and avoid sharp numerical gradients, we employed a smoothed transition based on the hyperbolic tangent function:

\[\displaystyle I_{\mathrm{stim}}(t) = 0.5\, I_{\mathrm{app}} \left[ 1 - \tanh\big(50\cdot(t -T_{\mathrm{sti}})\big) \right]. \]
This function provides a smooth approximation of a step-like stimulus while avoiding sharp temporal discontinuities that could affect numerical accuracy and solver stability.

\begin{figure}[t]
    \centering
    \includegraphics[width=\linewidth]{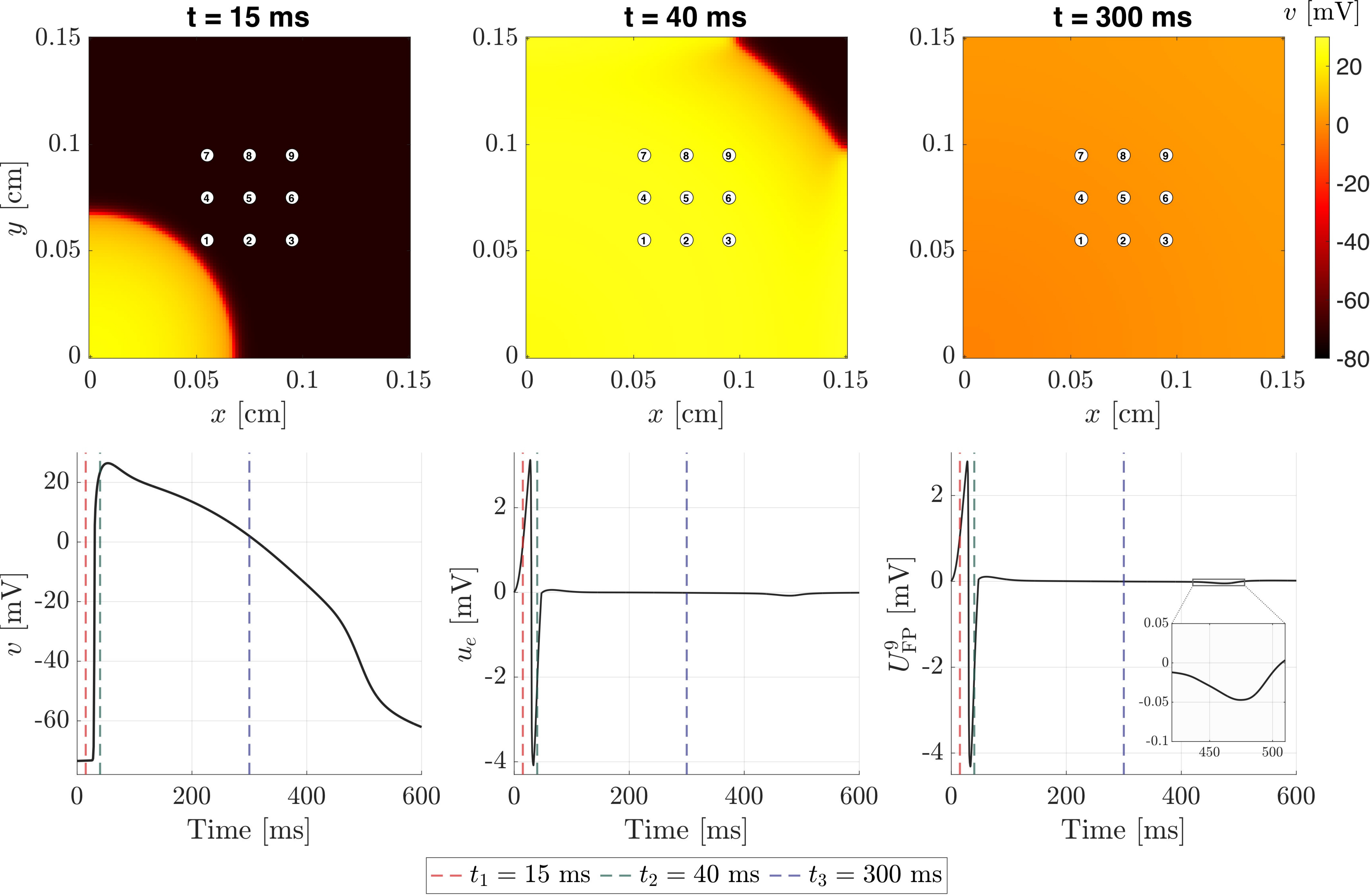}
    \caption{Space--time dynamics of the coupled MEA system. Top row: Spatial distribution of $v$ at representative time instants $t = 20,40,300\,\mathrm{ms}$. MEA electrodes are denoted by white markers and labelled according to their index. Bottom row: temporal evolution at electrode 9 of $v$ (left), $u_e$ (centre), and $U_{\mathrm{FP}}^9$ (right). Colour-coded vertical dashed lines denote the time points corresponding to the spatial snapshots. The inset in the $U_{\mathrm{FP}}$ plot provides a magnified view of the T-wave morphology.}
    \label{plot_v_FP}
\end{figure}

\section{Numerical Methods}
\label{sec:NumericalMethods}
To compute the numerical solution of the MEA system  (\ref{eq:mea_formulation}), we discretize the governing equations by bilinear finite elements in space, define a second-order Strang-based MEA operator splitting framework, and employ IMEX Runge-Kutta methods in time, described  in the following sections.

\subsection{MEA spatial discretization}\label{sp_disc}
The MEA system \eqref{eq:mea_formulation} is discretized in space using FEM with piecewise bilinear basis functions on a uniform squared mesh $\mathscr{T}^h$ of size $h = 15\,\mu\text{m}$ (\cite{Zienkiewicz_FEM}). Let $V_h = \text{span}\{\varphi_i\}_{i=1}^{N}$ be the corresponding finite element space. By applying a standard Galerkin procedure, we define the mass matrix $\mathrm{M}\in\mathbb{R}^{N\times N}$, the stiffness matrices $\mathrm{A}_{i,e}\in\mathbb{R}^{N\times N}$, and the electrode-specific local mass matrices $\mathrm{M}_k\in\mathbb{R}^{N\times N} \text{ for} \, k=1\ldots \Nel$:
    \[
    \mathrm{M}_{ij}=\int_\Omega \varphi_i\varphi_j\,dx,\quad
    (\mathrm{A}_{i,e})_{ij}=\int_\Omega \nabla\varphi_i\cdot\sigma_{i,e}\nabla\varphi_j\,dx,\quad
    (\mathrm{M}_k)_{ij}=\int_\Omega \varphi_i\varphi_j \chi_{e_k}\,dx.
    \]
The spatially semi-discrete problem consists in finding the vectors of nodal unknowns $\mathbf{v}(t)$ and $\mathbf{u}_e(t)$ such that:
\begin{equation} 
\label{eq:semi_discrete}
\begin{cases}
        \displaystyle
		C_m \mathrm{M}\frac{d\mathbf{v}}{dt} + \mathrm{A}_i \mathbf{v} + \mathrm{A}_i \mathbf{u}_e
        = I_{\mathrm{app}} - \mathrm{M}\,I_{\mathrm{ion}}(\mathbf{v},\mathbf{w},\mathbf{c}), \\
        \displaystyle
		(\mathrm{A}_i + \mathrm{A}_e)\mathbf{u}_e + \mathrm{A}_i \mathbf{v}
        = \frac{1}{\mathrm{z_{thick}}}\sum_{k=1}^{\Nel} \frac{1}{|e_k|}\mathrm{M}_k I_{el}^k, \qquad k=1\ldots \mathrm{N}_{el}
\end{cases}
\end{equation}
Together with the semi-discrete system \eqref{eq:semi_discrete}, also the average extracellular potential $U^k$ on each electrode, defined in \eqref{mea_3}, must be numerically computed. 
Specifically, at each time step, we first construct a scattered interpolant of the nodal values of $u_e$. This interpolant, denoted as $\hat{u}_e$, is then evaluated on a regular grid covering $e_k$. The integral is then computed using the composite 2D trapezoidal rule (see \cite{quarteroni}):
\[
U^k(t) \approx \frac{1}{|e_k|} \sum_{j=1}^{J_k} \omega_j \hat{u}_e(\boldsymbol{x}_j, t),
\]
where $\boldsymbol{x}_j$ are the $J_k$ integration points of the grid inside $e_k$ and $\omega_j$ are the corresponding quadrature weights. For a uniform grid with element area $h^2$, the trapezoidal weights are defined as $\omega_j = h^2$ for internal grid points, $\omega_j = h^2 / 2$ for points on the edges of $e_k$, and $\omega_j = h^2 / 4$ for the corner points.

\subsection{MEA Operator Splitting Framework}\label{OS}

To efficiently handle the different physical processes involved in the MEA formulation \eqref{eq:mea_formulation}, we adopt a second order operator splitting (OS) strategy (\cite{spiteri_second_order}). This fractional-step approach separates the integration of the reaction and diffusion operators, allowing each part to be treated with a dedicated numerical solver. Specifically, the time integration is performed through a second-order Strang-based OS scheme \cite{pavarino, sundnes}. Let $t_n$ and $t_{n+1} = t_n + \Delta t$ denote two consecutive time instants. Starting from the discrete solution at time $t_n$, the update at time $t_{n+1}$ is obtained through the following sequence of fractional steps.

\medskip
\textbf{Step 1 - Reaction half-step on $[t_n, t_n + \tfrac{\Delta t}{2}]$}\\
The ionic model is integrated over the first half of the time interval:
\begin{equation}
\label{step1}
\begin{cases}
C_m \mathrm{M} \dfrac{d\mathbf{v}}{dt} + \mathrm{M} I_{\mathrm{ion}}(\mathbf{v}, \mathbf{w}, \mathbf{c}) = 0, \\[8pt]
\dfrac{d\mathbf{w}}{dt} = \mathbf{W}(\mathbf{v}, \mathbf{w}, \mathbf{c}), \\
\dfrac{d\mathbf{c}}{dt} = \mathbf{C}(\mathbf{v}, \mathbf{w}, \mathbf{c}),
\end{cases}
\qquad t \in \left[t_n, t_n + \tfrac{\Delta t}{2}\right].
\end{equation}
Denoting by $\mathbf{v}^{*}, \mathbf{w}^{*}, \mathbf{c}^{*}$ the discrete quantities 
at time $t_n + \Delta t/2$, these intermediate values are used as initial data for the subsequent diffusion step.\\

\medskip
\noindent\textbf{Step 2 — Diffusion step on $[t_n, t_{n+1}]$}\\
Starting from $\mathbf{v}^{*}$, the pure diffusion problem is integrated over the full time interval:
\begin{equation}
\label{step2}
C_m \mathrm{M} \dfrac{d\mathbf{v}}{dt} + \mathrm{A}_i \mathbf{v} + \mathrm{A}_i \mathbf{u}_e = \mathrm{I}_{\text{app}}, \qquad t \in [t_n, t_{n+1}].
\end{equation}
The solution of this system provides the updated intermediate potential $\mathbf{v}^{**}$.\\

\medskip
\noindent\textbf{Step 3 — Reaction half-step on $[t_n + \tfrac{\Delta t}{2}, t_{n+1}]$}\\
The same ionic system of \hyperref[step1]{Step 1} is integrated again over the second half of the time interval. Starting from the intermediate values $\mathbf{v}^{**}, \mathbf{w}^{*}, \mathbf{c}^{*}$, this step provides the fully updated state variables:
\begin{equation*}
    \label{step3}
    \mathbf{v}(t_{n+1}), \quad \mathbf{w}(t_{n+1}), \quad \mathbf{c}(t_{n+1}).
\end{equation*}
\smallskip
\noindent\textbf{Step 4 — Extracellular and electrodes update}\\
Once $\mathbf{v}(t_{n+1})$ is known, $\mathbf{u}_e$ and $I_{el}^k$ must be updated by solving 
\begin{equation}
\label{step4}
\begin{cases}
\displaystyle
(\mathrm{A}_i + \mathrm{A}_e)\mathbf{u}_e + \mathrm{A}_i \mathbf{v}(t_{n+1}) = \displaystyle \frac{1}{z_{\text{thick}}} \sum_{k=1}^{\Nel} \frac{I_{el}^k}{|e_k|} \mathrm{M}_k, \\
\displaystyle \dfrac{d I_{el}^k}{dt} + \dfrac{I_{el}^k}{\tau} = \dfrac{C_{el}}{\tau} \dfrac{d U^k}{dt}, \qquad \text{with } U^k = \frac{1}{|e_k|} \int_{e_k} u_e(\boldsymbol{x}, t_{n+1}) \, de_k,
\end{cases}
\end{equation}
for $k = 1, \dots, \Nel$.

\subsection{MEA time discretization}

To advance the semi-discrete MEA system \eqref{eq:semi_discrete} in time, we employ a hybrid strategy designed to preserve the second-order accuracy of the splitting framework introduced in section \ref{OS}. Specifically, the reaction-diffusion processes (\hyperref[step1]{Steps 1--3}) are treated using an Implicit-Explicit Runge-Kutta (IMEX-RK) method, while the coupled electrode equations (\hyperref[step4]{Step 4}) are integrated using the Crank-Nicolson (CN) scheme.

\subsection*{IMEX--RK time discretization} The time integration of the reaction-diffusion system is performed using an IMEX--RK scheme. In this framework, the reaction and diffusion operators are treated with different numerical strategies: the nonlinear ionic dynamics (\hyperref[step1]{Step 1} and \hyperref[step3]{Step 3}) are integrated explicitly, while the linear diffusion operator (\hyperref[step2]{Step 2}) is handled implicitly.
Consider a generic system of the form
\[
\frac{dY}{dt} = F(Y) + G(Y),
\]
where $F(Y)$ represents the explicitly-treated contribution (i.e., \hyperref[step1]{Step 1} and \hyperref[step3]{Step 3}) and $G(Y)$ the implicitly-treated one (i.e., \hyperref[step2]{Step 2}). An $s$--stage IMEX Runge--Kutta scheme reads
\begin{subequations}
\begin{align}
\displaystyle
Y^{(i)} &= Y^n 
+ \Delta t \sum_{j=1}^{i-1} \tilde{a}_{ij} F\big(Y^{(j)}\big)
+ \Delta t \sum_{j=1}^{i} a_{ij} G\big(Y^{(j)}\big),
\quad i = 1,\dots,s
\label{eq:imex_stage}
\\
\displaystyle Y^{n+1} &= Y^n 
+ \Delta t \sum_{i=1}^{s} \tilde{b}_i F\big(Y^{(i)}\big)
+ \Delta t \sum_{i=1}^{s} b_i G\big(Y^{(i)}\big),
\label{eq:imex_update}
\end{align}
\end{subequations}
where the quantities $Y^{(i)}$, $i=1,\dots,s$, denote the approximating solution at time $t_n + c_i \Delta t$, while $Y^{n+1}$ the approximating solution at time $t_{n+1}$.

The method is uniquely defined by two $s\times s$ matrices: $\tilde{A}=(\tilde{a}_{ij})$, which is strictly lower triangular ($\tilde{a}_{ij}=0 \quad \text{for } j\ge i$) to explicitly treat $F$, and $A=(a_{ij})$, which defines the implicit discretization of $G$. The method is further characterized by the internal nodes
\[
\tilde{c}_i = \sum_{j=1}^{i-1} \tilde{a}_{ij}, 
\qquad 
c_i = \sum_{j=1}^{i} a_{ij},
\qquad i=1,\dots,s,
\]
and by the weight vectors $\tilde{b}, b \in \mathbb{R}^s$. These coefficients are conventionally summarized using the double Butcher tableau notation as reported in Table \ref{IMEX_Butcher}. Further details on IMEX-RK schemes can be found in \cite{Ascher_IMEX}, \cite{Boscarino_error} and \cite{Russo}, while for general RK schemes we refer to \cite{Hairer}.
\begin{table}[!b]
\small
\centering 
\caption{Double Butcher tableaus for IMEX-RK schemes. For each method, the left tableau ($\tilde{c}, \tilde{A}, \tilde{b}^{\top}$) defines the explicit discretization of the reaction term in \hyperref[step1]{Step 1} and \hyperref[step3]{Step 3}, while the right tableau ($c, A, b^{\top}$) defines the diagonally implicit discretization of the diffusion operator in \hyperref[step2]{Step 2}.}
\label{IMEX_Butcher}
\vspace{5pt}
\renewcommand{\arraystretch}{1.5} 

\begin{tabular}{cc}
\textbf{\small General IMEX-RK Form} & \textbf{\small H(2,2,2)} \\
\addlinespace[5pt]
$\begin{array}{c|c} \tilde{c} & \tilde{A} \\ \hline & \tilde{b}^{\top} \end{array} \quad \begin{array}{c|c} c & A \\ \hline & b^{\top} \end{array}$ 
& 
$\begin{array}{c|cc} 0 & 0 & 0 \\ 1 & 1 & 0 \\ \hline & 1/2 & 1/2 \end{array} \quad 
\begin{array}{c|cc} \frac{2-\sqrt{2}}{2} & \frac{2-\sqrt{2}}{2} & 0 \\ \frac{\sqrt{2}}{2} & \sqrt{2}-1 & \frac{2-\sqrt{2}}{2} \\ \hline & 1/2 & 1/2 \end{array}$ \\

\addlinespace[15pt]

& \textbf{\small SSP2(2,2,2)} \\
\addlinespace[5pt]
& $\begin{array}{c|cc} 0 & 0 & 0 \\ 1 & 1 & 0 \\ \hline & 1/2 & 1/2 \end{array} \quad 
\begin{array}{c|cc} 1-\frac{\sqrt{2}}{2} & 1-\frac{\sqrt{2}}{2} & 0 \\ 1 & \frac{\sqrt{2}}{2} & 1-\frac{\sqrt{2}}{2} \\ \hline & \frac{\sqrt{2}}{2} & 1-\frac{\sqrt{2}}{2} \end{array}$
\end{tabular}

\end{table}

\subsection*{Crank-Nicolson method for electrode dynamics}
Once the updated state variables from the reaction-diffusion integration (\hyperref[step1]{Steps 1--3}) are computed, the coupled electrode-extracellular system must be advanced (\hyperref[step4]{Step 4}). 
To maintain consistency with the overall second-order accuracy of the OS strategy and remaining fully consistent with the implicit treatment of the coupled system \eqref{step4} we employ the CN scheme. 
Over the interval $[t_n,t_{n+1}]$, the discretized system reads:
\begin{equation}
\label{eq:CN_electrode}
\begin{dcases}
        \left(\mathrm{A}_i+\mathrm{A}_e\right)\mathbf{u_e}^{n+1}-\frac{1}{z_{\text{thick}}}\sum_{k=1}^{\Nel}\frac{I_{el}^{k,n+1}}{|e_k|}\mathrm{M}_k=-\mathrm{A}_i\mathbf{v}^{n+1}\\[6pt]
        \left(\frac{1}{\Delta t}+\frac{1}{2\tau}\right)I_{el}^{k,n+1}-\frac{C_{el}}{\tau\Delta t}U^{k,n+1} = \left(\frac{1}{\Delta t}-\frac{1}{2\tau}\right) I_{el}^{k,n} -\frac{C_{el}}{\tau\Delta t}U^{k,n}.
\end{dcases}
\end{equation}
for $k = 1, \dots, \Nel$. The resulting problem is therefore solved as a coupled linear system for the unknowns $\mathbf{u_e}^{n+1}$ and $I_{el}^{k,n+1}$ for $k=1,\dots,\Nel$.

\subsection{Error computation and order of convergence}
\label{sec:error_analysis}
To verify the numerical implementation and assess the performance of the proposed framework, we conduct a temporal convergence study. In this section we describe the reference solution, the error metrics and the procedure used to estimate the experimental order of convergence.

\begin{figure}[tbp]
    \centering
    \includegraphics[width=\textwidth]{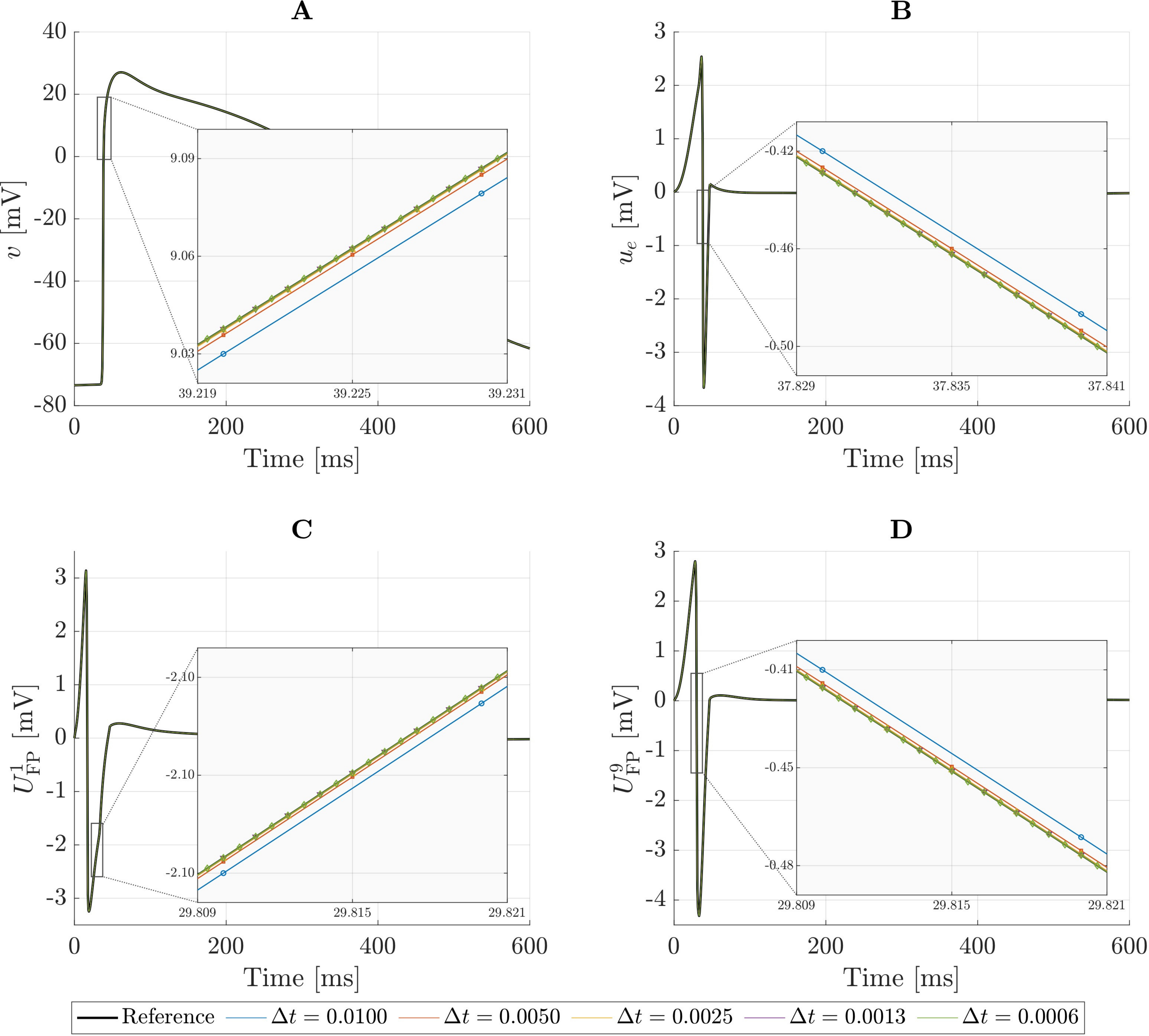}
    \caption{Temporal convergence analysis of the simulated electrophysiological signals for the 2D MEA model using the H(2,2,2) numerical scheme. The panels display $v$ (A) and $u_e$ (B) computed at a point located at 80\% of the domain length, and the field potentials $U_{\mathrm{FP}}$ at electrode 1 (C) and electrode 9 (D). In each plot, the solid black line represents the reference solution, while the coloured lines correspond to the various $\Delta t$ tested. The magnified insets demonstrate the consistent convergence of the numerical results toward the reference solution as the temporal resolution is increased.}
    \label{fig:convergence_study}
\end{figure}

\subsection*{Reference solution} Since analytical solutions for the coupled MEA system are not available, a reference solution $(v_{\mathrm{ref}}, u_{e,\mathrm{ref}}, I_{el,\mathrm{ref}}^k)$ is computed using the second-order IMEX-SSP2(2,2,2) scheme (see Table \ref{IMEX_Butcher}). This scheme is selected for its high stability and accuracy in treating stiff reaction--diffusion systems (see \cite{Russo} for more details).
To ensure that the reference solution is sufficiently accurate, we use a very fine time step:

\[
\Delta t_{\mathrm{ref}} =
\frac{10^{-2}}{2^{6}}
=
1.5625 \cdot 10^{-4}\ \mathrm{ms},
\]
This choice makes the temporal discretization error of the reference solution negligible compared to the ones measured in the following convergence study.

\subsection*{Accuracy evaluation}
To assess the accuracy of the proposed method, two different error measures were considered. First, we employed the mixed root mean square (MRMS) error introduced in \cite{cervi2018high}. 

Let's consider the fully discrete space--time set obtained by combining all the nodes of the spatial grid and all the discrete time instants of the simulation. Let $N = N_{\mathrm{S}} \times N_{\mathrm{T}}$ be the total number of these space--time points, where $N_{\mathrm{S}}$ and $N_{\mathrm{T}}$ denote the number of spatial grid nodes and discrete time instants, respectively. Let now denote with $q_i$ and $q_{\mathrm{ref},i}$ the MEA system quantities (i.e.\ $v$, $u_e$ or \FP[k] for $k=1,\ldots,\mathrm{N}_{\mathrm{el}}$) computed by the test and reference solutions, respectively, and evaluated at the $i$-th space--time point. 
The MRMS error is then defined as:

\begin{equation}
\label{eq:mrms}
\MRMS_q =\sqrt{\frac{1}{N}\sum_{i=1}^{N}\left(\frac{q_{\mathrm{ref}_i} - q_i}{1 + |q_{\mathrm{ref}}|}\right)^2}.
\end{equation}

In addition to this global measure, we also evaluated the error using discrete $L^2$ norm computed at a fixed spatial location and a fixed time instant. More precisely, for a selected spatial node $x^*$ of the computational grid, the discrete temporal $L^2$ error for the variable $q$ is:
\[
\displaystyle
\err[q](x^*) =\sqrt{\frac{1}{N_t}\sum_{n=1}^{N_t}\big(q(x^*,t_n) - q_{\mathrm{ref}}(x^*,t_n)\big)^2},
\label{space_err}
\]
While for a selected time instant $t^*$, the spatial $L^2$ error for the variable $q$ is:
\[
\displaystyle
\err[q](t^*) = \sqrt{\frac{1}{N_s}\sum_{n=1}^{N_s}\big(q(x_n,t^*) - q_{\mathrm{ref}}(x_n,t^*)\big)^2}.
\]

\subsection*{Order of Convergence}
The global order of convergence \ooc{real} of the framework depends on both the internal order of the IMEX--RK method ($\ooc{\mathrm{IMEX}}$) and the order of the OS procedure ($\ooc{\mathrm{split}}$). Thus, the theoretical overall temporal order satisfies:
\[
\ooc{real} = \min(\ooc{IMEX},\, \ooc{split}).
\]
Since our OS splitting is based on the Strang scheme, which is second-order accurate in time (see \cite{sundnes}), whenever $\ooc{IMEX} \ge 2$ the combined method is expected to achieve second-order temporal accuracy. 

To verify this, the experimental order of convergence \ooc{} is estimated by comparing errors across a sequence of refined time steps. Starting from an initial $\Delta t$, a sequence of numerical solutions was computed by repeatedly halving the time step, providing $\Delta t_k = \Delta t / 2^{k-1}$ for $k=1, 2, \dots$. For each refinement level $k$, the error \err[k] is computed using one of the strategies introduced in the previous paragraph. The experimental order of convergence between two consecutive levels is then given by:

\begin{equation}
\label{eq:eoc}
\displaystyle
\ooc{\err[q]} =\frac{\log(\err[q,1]/\err[q,2])}{\log \left(\Delta t_1/\Delta t_2\right)}.
\end{equation}

Where \err[q,1] and \err[q,2] denote the errors associated with the time steps $\Delta t_1$ and $\Delta t_2$ respectively and computed using both the MRMS error and the $L^2$ errors for $v$, $u_e$, and $\FP[5]$, ensuring a comprehensive assessment of the convergence.

\subsection*{Implementation details} 
All the simulations are conducted using the high-performance computing cluster of the University of Pavia (FAT nodes). Each simulation was executed on a single compute node equipped with 32 physical cores and 256 GB of RAM. Parallelization was achieved via shared-memory within the \textsc{Matlab} environment. The simulated time interval for each run was set to $600\,\mathrm{ms}$, providing a sufficiently long window to observe a full AP dynamic of a ventricular--like hiPSC-CM. Table \ref{tab:model_parameters} summarizes all the parameters adopted in the simulations.

\begin{table}[b]
\renewcommand{\arraystretch}{1.2}
\centering
\caption{Physical parameters, model constants and discretization settings used in the simulations.}
\label{tab:model_parameters}
\begin{tabular}{lll}
\toprule
\textbf{Parameter} & \textbf{Symbol} & \textbf{Value} \\ 
\midrule
Membrane capacitance \cite{pavarino} & $C_m$ & $1.0\ \mu\mathrm{F\ cm^{-2}}$ \\
Intracellular conductivity \cite{Tveito_2} & $\sigma_i$ & $0.01\ \mathrm{mS\ cm^{-1}}$ \\
Extracellular conductivity \cite{Tveito_2} & $\sigma_e$ & $0.1\ \mathrm{mS\ cm^{-1}}$  \\
Electrode capacitance & $C_{el}$ & $10^{-10}\ \mathrm{F}$  \\
Electrode ground resistance \cite{Moulin} & $R_i$ & $10^{9}\ \Omega$ \\
Electrode internal resistance \cite{Moulin} & $R_{el}$ & $10^{6}\ \Omega$ \\
Extracellular layer thickness \cite{Coudier} & $\mathrm{z}_{\text{thick}}$ & $1\ \mu\mathrm{m}$ \\
Spatial mesh size & $h$ & $15\ \mu\mathrm{m}$ \\ 
Applied current density & $I_{\mathrm{app}}$ & $150\ \mu\mathrm{A\ cm^{-2}}$ \\
\bottomrule
\end{tabular}
\end{table}

\section{Numerical Results and accuracy study}
\label{sec:results}

\subsection{0D Paci 2020 Ionic Model}
\label{sec:0d_paci}

Before presenting the full 2D tissue simulations, we assess the numerical performance, temporal accuracy, and computational efficiency of the time-stepping schemes using a 0D single-cell configuration. 

We test the explicit parts of the three numerical schemes used in the 2D configuration, namely the $\mathrm{SP}(1,1,1)$ scheme, acting as an Explicit Euler method (EE), and the $\mathrm{H}(2,2,2)$ and $\mathrm{SSP2}(2,2,2)$ schemes, which share the same explicit Butcher tableau of Heun's method. Since their convergence behaviours are identical in this 0D setup, their performance is reported under a single category, denoted as Heun. All convergence measurements are computed against a refined reference solution generated using Heun's method with:
$$
\Delta t_{\mathrm{ref}} = \frac{10^{-2}}{2^{6}} = 1.5625 \cdot 10^{-4}\ \mathrm{ms}.
$$

\begin{figure}[tbp]
    \centering
    \includegraphics[width=\linewidth]{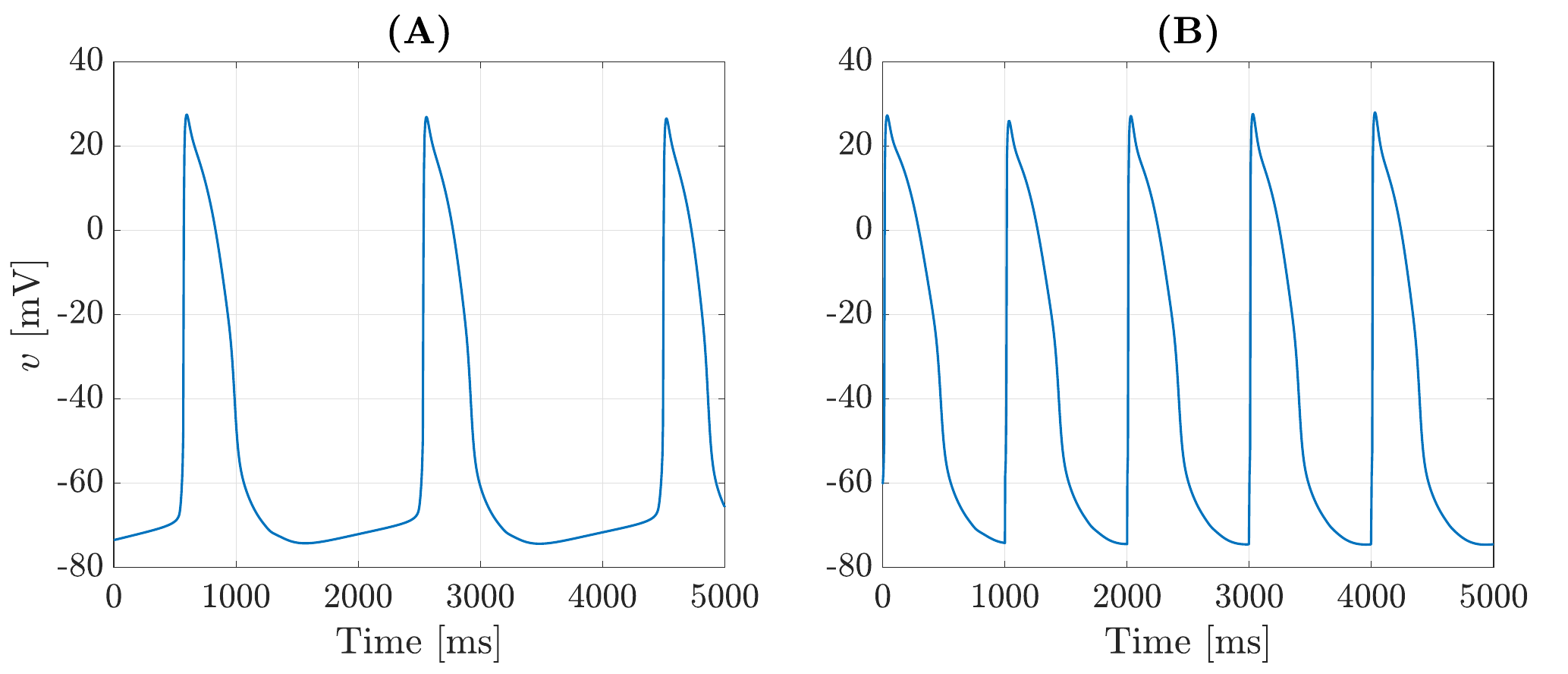}
    \caption{Temporal evolution of $v$ in a single hiPSC-CM using Paci 2020 ionic model under unpaced conditions (A) and within $1\ \mathrm{Hz}$ electrical pacing (B).}
    \label{plot_v_Cai}
\end{figure}

The simulations are carried out for a total of $5$ s under two distinct frameworks: the unpaced and the paced one. In the unpaced configuration, no stimulus is given and the Paci 2020 model simulates the typical spontaneous electrical activity of hiPSC-CMs. Conversely, the paced configuration implements a forced pacing strategy at a regular frequency of $1\ \mathrm{Hz}$. Pacing was represented as an instantaneous depolarization of the membrane potential at the prescribed stimulation times. This approach can be interpreted as an infinitely short stimulus pulse which bring the membrane potential above the excitation threshold, developing an AP according to the intrinsic ionic dynamics of the model \cite{Vinet1994}. Figure \ref{plot_v_Cai} displays the regular train of APs for both configurations. The initial conditions used are reported in \cite{Paci2020}.

\subsection*{Error and Order of Convergence}
\label{Convergence_0d}

We investigate the temporal accuracy of the proposed numerical schemes by computing the discrete temporal $L^2$ error for $v$ and the intracellular calcium concentration ($Ca_i$). The errors are computed across a sequence of systematically halved time steps starting from $\Delta t = 10^{-2} \mathrm{ms}$. The numerical errors and the corresponding experimental orders of convergence $\ooc{}$ are reported in Table \ref{tab:err_0d}a for the unpaced configuration, and in Table \ref{tab:err_0d}b for the paced one. As expected, the numerical error decreases across all tested methods when $\Delta t$ is reduced, exhibiting the correct theoretical order of convergence in both the unpaced and paced framework.

\begin{table}[htbp]
    \centering
    \setlength{\tabcolsep}{3pt}
    \caption{Temporal $L^2$ errors and experimental orders of convergence ($\ooc{}$) for $v$ and $Ca_i$ in the 0D hiPSC-CM ionic model, computed for different time steps $\Delta t$. Results are shown for the unpaced (a) and paced (b) configurations.}
    
    \textbf{(a) Unpaced configuration} \\[6pt] 
    \begin{tabular}{l cccc cccc}
    \toprule
    \multirow{3}{*}{$\Delta t$ [ms]} & \multicolumn{4}{c}{EE} & \multicolumn{4}{c}{Heun} \\
    \cmidrule(lr){2-5} \cmidrule(lr){6-9}
     & $\err[v]$ & $\ooc{v}$ & $\err[Ca_i]$ & $\ooc{Ca_i}$ & $\err[v]$ & $\ooc{v}$ & $\err[Ca_i]$ & $\ooc{Ca_i}$ \\
    \midrule
    1e-2    & 2.39e-2 & --    & 3.19e-6 & --    & 3.24e-5 & --    & 3.70e-9 & --    \\
    5e-3    & 1.19e-2 & 1.008 & 1.59e-6 & 1.003 & 8.07e-6 & 2.003 & 1.00e-9 & 1.885 \\
    2.5e-3  & 5.96e-3 & 1.004 & 7.95e-7 & 1.001 & 2.01e-6 & 2.005 & 2.59e-10& 1.950 \\
    1.25e-3 & 2.98e-3 & 1.002 & 3.97e-7 & 1.001 & 4.97e-7 & 2.017 & 6.53e-11& 1.990 \\
    6.25e-4 & 1.49e-3 & 1.001 & 1.99e-7 & 1.001 & 1.18e-7 & 2.071 & 1.57e-11& 2.055 \\
    \bottomrule
    \end{tabular}
    
    \vspace{20pt} 
    
    \textbf{(b) Paced configuration} \\[6pt] 
    \begin{tabular}{l cccc cccc}
    \toprule
    \multirow{3}{*}{$\Delta t$ [ms]} & \multicolumn{4}{c}{EE} & \multicolumn{4}{c}{Heun} \\
    \cmidrule(lr){2-5} \cmidrule(lr){6-9}
     & $\err[v]$ & $\ooc{v}$ & $\err[Ca_i]$ & $\ooc{Ca_i}$ & $\err[v]$ & $\ooc{v}$ & $\err[Ca_i]$ & $\ooc{Ca_i}$ \\
    \midrule
    1e-2    & 2.68e-2 & --    & 1.19e-5 & --    & 3.16e-4 & --    & 1.75e-8 & --    \\
    5e-3    & 1.36e-2 & 0.980 & 5.89e-6 & 1.010 & 8.98e-5 & 1.917 & 5.55e-9 & 1.651 \\
    2.5e-3  & 6.85e-3 & 0.990 & 2.94e-6 & 1.005 & 2.36e-5 & 1.931 & 1.52e-9 & 1.866 \\
    1.25e-3 & 3.44e-3 & 0.995 & 1.47e-6 & 1.002 & 5.95e-6 & 1.984 & 3.93e-10& 1.956 \\
    6.25e-4 & 1.72e-3 & 0.998 & 7.32e-7 & 1.001 & 1.43e-6 & 2.055 & 9.55e-11& 2.041 \\
    \bottomrule
    \end{tabular}

    \vspace{12pt}
    
    \label{tab:err_0d} 
\end{table}

\subsection*{Comparative Performance Analysis}
\label{subsec:performance0d}

To complement the convergence study, we evaluate the overall computational efficiency of the numerical methods by studying the correlation between numerical error, CPU time, and number of time steps for $v$ and $Ca_i$. The CPU times measured for each solver are summarized in Table \ref{tab:cpu_time_0d}a for the unpaced framework and in Table \ref{tab:cpu_time_0d}b for the paced one. 

\begin{figure}[htbp]
    \centering
    \includegraphics[width=\linewidth]{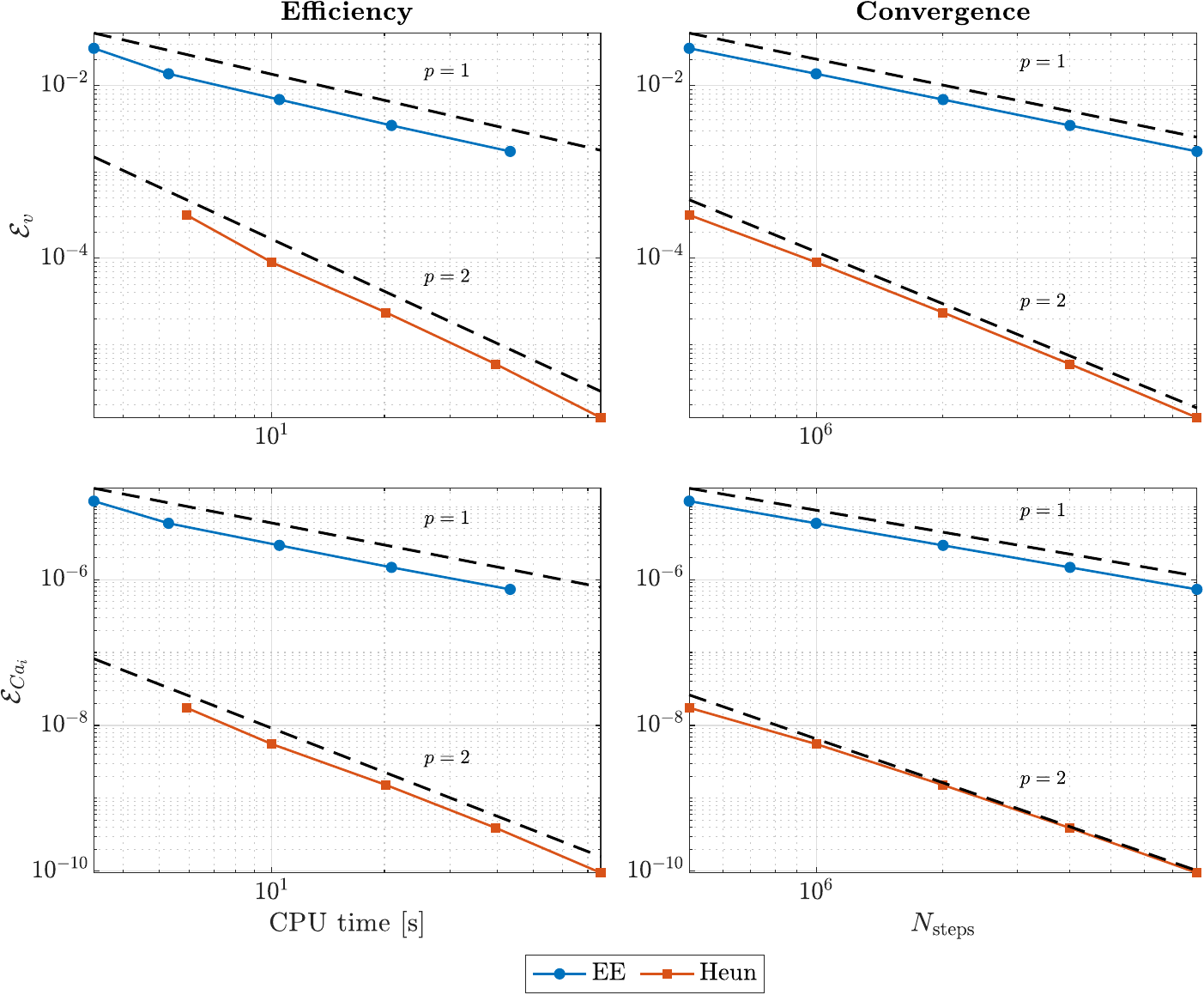}
    \caption{Performance and convergence analysis for the 0D Paci 2020 ionic model under paced conditions. The left column shows the computational efficiency evaluated as $\err[]$ vs. CPU time. The right column shows the asymptotic convergence profile as error vs. number of time steps ($N_{\mathrm{steps}}$). Top row refer to $v$, bottom row to $Ca_i$. Solid lines with markers show the experimental data, while dashed lines represent the expected theoretical slopes ($p=1$ and $p=2$). Results for the unpaced configuration are entirely comparable to these ones.}
    \label{fig:Comp_performance_0d}
\end{figure}

Figure \ref{fig:Comp_performance_0d} presents the results for the paced configuration on log-log scales to emphasize the asymptotic behaviour of each solver and facilitate a direct comparison. The results for the unpaced case are entirely comparable to the paced one and, thus, not explicitly reported.

\begin{table}[htbp]
    \centering
    \small
    
    \caption{Computational CPU time and numerical accuracy $\err[v]$ for decreasing values of $\Delta t$ in the 0D hiPSC-CM ionic model. Results a shown for the unpaced (a) and paced (b) configurations.}
    \label{tab:cpu_time_0d}
    \setlength{\tabcolsep}{2pt}
    
    \begin{tabular}[t]{c}
        \textbf{(a) Unpaced configuration} \\[6pt]
        \begin{tabular}{l cc cc}
        \toprule
        \multirow{2}{*}{$\Delta t$ [ms]} & \multicolumn{2}{c}{EE} & \multicolumn{2}{c}{Heun} \\
        \cmidrule(lr){2-3} \cmidrule(lr){4-5}
         & CPU [s] & $\err[v]$ & CPU [s] & $\err[v]$ \\
        \midrule
        1e-2    & 4.05  & 2.39e-2 & 6.35  & 3.24e-5 \\
        5e-3    & 6.34  & 1.19e-2 & 10.25 & 8.07e-6 \\
        2.5e-3  & 11.67 & 5.96e-3 & 21.56 & 2.01e-6 \\
        1.25e-3 & 22.15 & 2.98e-3 & 37.46 & 4.97e-7 \\
        6.25e-4 & 44.75 & 1.49e-3 & 74.98 & 1.18e-7 \\
        \bottomrule
        \end{tabular}
    \end{tabular}%
    \qquad 
    \begin{tabular}[t]{c}
        \textbf{(b) Paced configuration} \\[6pt]
        \begin{tabular}{l cc cc}
        \toprule
        \multirow{2}{*}{$\Delta t$ [ms]} & \multicolumn{2}{c}{EE} & \multicolumn{2}{c}{Heun} \\
        \cmidrule(lr){2-3} \cmidrule(lr){4-5}
         & CPU [s] & $\err[v]$ & CPU [s] & $\err[v]$ \\
        \midrule
        1e-2    & 1.12  & 2.68e-2 & 4.37  & 3.16e-4 \\
        5e-3    & 2.45  & 1.36e-2 & 4.81  & 8.98e-5 \\
        2.5e-3  & 4.90  & 6.85e-3 & 8.84  & 2.36e-5 \\
        1.25e-3 & 9.85  & 3.44e-3 & 19.62 & 5.95e-6 \\
        6.25e-4 & 19.80 & 1.72e-3 & 38.73 & 1.43e-6 \\
        \bottomrule
        \end{tabular}
    \end{tabular}
\end{table}

\subsection{2D MEA results}
In this section, we assess the numerical performance of the 2D MEA model considering only the paced configuration. The unpaced case is not considered since the cells would rapidly synchronize their activity, causing the entire domain to depolarize almost simultaneously and preventing the formation of propagation wavefronts. Furthermore, in order to limit the total computational cost of the simulations, the analysis is restricted to a single paced AP ($\approx 600,\mathrm{ms}$).

\subsection*{APs and space--time propagation}

To illustrate the electrophysiological dynamics captured by the proposed framework, Fig.~\ref{plot_v_FP} provides an overview of the spatio-temporal evolution of the main model variables. Following the application of the stimulus, a depolarization wavefront propagates smoothly across the tissue domain. This activation process is clearly observed in the spatial distribution of the transmembrane potential $v$ (top row of Fig.~\ref{plot_v_FP}), shown at three representative time instants.

The bottom row of Fig.~\ref{plot_v_FP} reports the temporal traces of $v$, $u_e$, and \FP[9] recorded at the ninth electrode. The highlighted markers establish the correspondence between the temporal signals and the snapshots displayed in the top row. All quantities exhibit the expected qualitative behaviour, demonstrating the ability of the model to reproduce both the propagation dynamics within the tissue and the extracellular signals measured by the MEA.

\subsection*{Error and order of convergence}

We investigate the temporal accuracy of the proposed numerical schemes by comparing the results with the reference solution described in Section \ref{sec:error_analysis}. The analysis is performed for $v$, $u_e$, and $U_{\mathrm{FP}}^5$. In particular, for the discrete $L^2$ norms, we consider the centre node of the domain for the error in time $\err[1](x^*)$ and the time instant $t=300$~ms for the error in space $\err[q](t^*)$.

The error measured in terms of \MRMS, $\err[q](x^*)$, and $\err[q](t^*)$ for decreasing values of $\Delta t$ are reported in Table \ref{tab:main_errors}. For \FP[5], the error is evaluated only in time, since this quantity is defined at discrete electrode locations and does not admit a continuous spatial norm. The corresponding experimental orders of convergence, computed according to \eqref{eq:eoc}, are collected in Tables \ref{tab:main_orders}\, (a)-(c). We compare three IMEX--RK schemes: the SP(1,1,1) scheme, which combines Forward Euler for the explicit part and Backward Euler for the implicit part, the IMEX--H(2,2,2) scheme, coupling a second-order explicit RK discretization with a CN scheme, and the SSP2(2,2,2) scheme, already described in Section \ref{sec:error_analysis}. For the specific Butcher tableaus of the proposed methods see \ref{IMEX_Butcher}.

The results show a consistent reduction of the error as the time step decreases for all quantities and numerical methods. This convergence behaviour is qualitatively illustrated in Figure \ref{fig:convergence_study}, where the numerical solutions obtained with the IMEX--H(2,2,2) scheme are shown to progressively overlap with the reference traces as $\Delta t$ is refined. In particular, as expected from theory, the SP(1,1,1) scheme exhibits a first-order convergence behaviour, while both IMEX--H(2,2,2) and SSP2(2,2,2) display a second-order convergence trend across all error measures.

\begin{table}[tbp]
\centering
\setlength{\tabcolsep}{3pt}
\caption{Error measures of $v$ (a), $u_e$ (b) and \FP[5] (c) for the 2D MEA model. For each quantity, errors are reported for different IMEX-RK schemes as a function of $\Delta t$. For $v$ and $u_e$ errors are computed using MRMS, $\err[v](x^*)$, where $x^*$ is the central node of the domain, and $\err[v](t^*)$ with $t^*=300 \ \mathrm{ms}$. For $\FP[5]$ the error is evaluated using $\err[\FP](x^*)$.}
\label{tab:main_errors}

{\small \textbf{(a)} Error measures for $v$}\\[5pt]
\begin{tabular}{l ccc ccc ccc}
\toprule
\multirow{2}{*}{$\Delta t$ [ms]} & \multicolumn{3}{c}{SP(1,1,1)} & \multicolumn{3}{c}{H(2,2,2)} & \multicolumn{3}{c}{SSP2(2,2,2)} \\
\cmidrule(lr){2-4} \cmidrule(lr){5-7} \cmidrule(lr){8-10}
& \MRMS$_v$ & $\err[v](x^*)$ & $\err[v](t^*)$ & \MRMS$_v$ & $\err[v](x^*)$ & $\err[v](t^*)$ & \MRMS$_v$ & $\err[v](x^*)$ & $\err[v](t^*)$ \\
\midrule
1e-2 & 5.22e-1 & 4.37e-1 & 2.68e-2 & 2.59e-3 & 2.16e-3 & 1.23e-4 & 7.63e-3 & 6.72e-3 & 2.90e-4 \\
5e-3 & 2.65e-1 & 2.20e-1 & 1.34e-2 & 6.57e-4 & 5.47e-4 & 3.08e-5 & 1.93e-3 & 1.70e-3 & 7.28e-5 \\
2.5e-3 & 1.33e-1 & 1.10e-1 & 6.71e-3 & 1.64e-4 & 1.36e-4 & 7.65e-6 & 4.83e-4 & 4.25e-4 & 1.82e-5 \\
1.25e-3 & 6.66e-2 & 5.53e-2 & 3.36e-3 & 3.98e-5 & 3.30e-5 & 1.86e-6 & 1.20e-4 & 1.05e-4 & 4.50e-6 \\
6.25e-4 & 3.33e-2 & 2.77e-2 & 1.68e-3 & 8.55e-6 & 7.01e-6 & 4.12e-7 & 2.85e-5 & 2.51e-5 & 1.07e-6 \\
\bottomrule
\end{tabular}

\vspace{20pt} 

{\small \textbf{(b)} Error measures for $u_e$}\\[5pt]
\begin{tabular}{l ccc ccc ccc}
\toprule
\multirow{2}{*}{$\Delta t$ [ms]} & \multicolumn{3}{c}{SP(1,1,1)} & \multicolumn{3}{c}{H(2,2,2)} & \multicolumn{3}{c}{SSP2(2,2,2)} \\
\cmidrule(lr){2-4} \cmidrule(lr){5-7} \cmidrule(lr){8-10}
& \MRMS$_{u_e}$ & $\err[u_e](x^*)$ & $\err[u_e](t^*)$ & \MRMS$_{u_e}$ & $\err[u_e](x^*)$ & $\err[u_e](t^*)$ & \MRMS$_{u_e}$ & $\err[u_e](x^*)$ & $\err[u_e](t^*)$ \\
\midrule
1e-2 & 4.22e-2 & 3.89e-2 & 1.91e-4 & 2.08e-4 & 1.91e-4 & 8.85e-7 & 6.16e-4 & 5.96e-4 & 2.79e-6 \\
5e-3 & 2.14e-2 & 1.96e-2 & 9.60e-5 & 5.28e-5 & 4.85e-5 & 2.25e-7 & 1.55e-4 & 1.50e-4 & 7.05e-7 \\
2.5e-3 & 1.08e-2 & 9.85e-3 & 4.81e-5 & 1.32e-5 & 1.21e-5 & 5.62e-8 & 3.89e-5 & 3.77e-5 & 1.77e-7 \\
1.25e-3 & 5.39e-3 & 4.93e-3 & 2.41e-5 & 3.19e-6 & 2.92e-6 & 1.36e-8 & 9.64e-6 & 9.33e-6 & 4.37e-8 \\
6.25e-4 & 2.70e-3 & 2.47e-3 & 1.21e-5 & 6.85e-7 & 6.21e-7 & 2.89e-9 & 2.30e-6 & 2.22e-6 & 1.04e-8 \\
\bottomrule
\end{tabular}

\vspace{20pt}

{\small \textbf{(c)} Error measures for \FP[5]}\\[5pt]
\begin{tabular}{l ccc}
\toprule
\multirow{2}{*}{$\Delta t$ [ms]} & \multicolumn{3}{c}{$\err[\FP[5]](x^*)$} \\
\cmidrule(lr){2-4}
& SP(1,1,1) & H(2,2,2) & SSP2(2,2,2) \\
\midrule
1e-2 & 2.72e-2 & 1.34e-4 & 4.19e-4 \\
5e-3 & 1.37e-2 & 3.43e-5 & 1.06e-4 \\
2.5e-3 & 6.87e-3 & 8.56e-6 & 2.65e-5 \\
1.25e-3 & 3.44e-3 & 2.07e-6 & 6.56e-6 \\
6.25e-4 & 1.72e-3 & 4.40e-7 & 1.56e-6 \\
\bottomrule
\end{tabular}
\end{table}

\begin{table}[tbp]
\centering
\setlength{\tabcolsep}{3pt}
\caption{Experimental order of convergence $p$ for $v$ (a), $u_e$ (b) and $\FP[5]$ (c) in the 2D MEA model. The orders are calculated using the MRMS error ($\ooc{MRMS}$), the temporal error ($\ooc{\err[q](x^*)}$), and the spatial error ($\ooc{\err[q](t^*)}$).}
\label{tab:main_orders}

{\small \textbf{(a)} Experimental order of convergence for $v$}\\[5pt]
\begin{tabular}{l ccc ccc ccc}
\toprule
\multirow{2}{*}{$\Delta t$ [ms]} & \multicolumn{3}{c}{SP(1,1,1)} & \multicolumn{3}{c}{H(2,2,2)} & \multicolumn{3}{c}{SSP2(2,2,2)} \\
\cmidrule(lr){2-4} \cmidrule(lr){5-7} \cmidrule(lr){8-10}
& $\ooc{MRMS}$ & $\ooc{\err[v](x^*)}$ & $\ooc{\err[v](t^*)}$ & $\ooc{MRMS}$ & $\ooc{\err[v](x^*)}$ & $\ooc{\err[v](t^*)}$ & $\ooc{MRMS}$ & $\ooc{\err[v](x^*)}$ & $\ooc{\err[v](t^*)}$ \\
\midrule
1e-2    & - & - & - & - & - & - & - & - & - \\
5e-3    & 0.980 & 0.988 & 0.998 & 1.979 & 1.981 & 1.999 & 1.986 & 1.986 & 1.994 \\
2.5e-3  & 0.993 & 0.996 & 0.999 & 2.001 & 2.003 & 2.009 & 1.997 & 1.997 & 2.001 \\
1.25e-3 & 0.997 & 0.998 & 0.999 & 2.045 & 2.049 & 2.040 & 2.013 & 2.013 & 2.015 \\
6.25e-4 & 0.999 & 0.999 & 0.999 & 2.219 & 2.235 & 2.174 & 2.069 & 2.069 & 2.069 \\
\bottomrule
\end{tabular}

\vspace{20pt}

{\small \textbf{(b)} Experimental order of convergence for $u_e$}\\[5pt]
\begin{tabular}{l ccc ccc ccc}
\toprule
\multirow{2}{*}{$\Delta t$ [ms]} & \multicolumn{3}{c}{SP(1,1,1)} & \multicolumn{3}{c}{H(2,2,2)} & \multicolumn{3}{c}{SSP2(2,2,2)} \\
\cmidrule(lr){2-4} \cmidrule(lr){5-7} \cmidrule(lr){8-10}
& $\ooc{MRMS}$ & $\ooc{\err[u_e](x^*)}$ & $\ooc{\err[u_e](t^*)}$ & $\ooc{MRMS}$ & $\ooc{\err[u_e](x^*)}$ & $\ooc{\err[u_e](t^*)}$ & $\ooc{MRMS}$ & $\ooc{\err[u_e](x^*)}$ & $\ooc{\err[u_e](t^*)}$ \\
\midrule
1e-2    & - & - & - & - & - & - & - & - & - \\
5e-3    & 0.979 & 0.988 & 0.992 & 1.980 & 1.981 & 1.977 & 1.986 & 1.986 & 1.986 \\
2.5e-3  & 0.993 & 0.996 & 0.996 & 2.001 & 2.003 & 1.999 & 1.997 & 1.997 & 1.997 \\
1.25e-3 & 0.997 & 0.998 & 0.998 & 2.045 & 2.049 & 2.047 & 2.013 & 2.013 & 2.013 \\
6.25e-4 & 0.999 & 0.999 & 0.999 & 2.221 & 2.235 & 2.234 & 2.069 & 2.069 & 2.069 \\
\bottomrule
\end{tabular}

\vspace{20pt}

{\small \textbf{(c)} Experimental order of convergence for \FP[5]}\\[5pt]
\begin{tabular}{l ccc}
\toprule
\multirow{2}{*}{$\Delta t$ [ms]} & \multicolumn{3}{c}{$\ooc{\err[\FP[5]](x^*)}$} \\
\cmidrule(lr){2-4}
& SP(1,1,1) & H(2,2,2) & SSP2(2,2,2) \\
\midrule
1e-2    & - & - & -  \\
5e-3    & 0.992 & 1.982 & 1.986 \\
2.5e-3  & 0.996 & 2.003 & 1.997 \\
1.25e-3 & 0.998 & 2.048 & 2.014 \\
6.25e-4 & 0.999 & 2.233 & 2.069 \\
\bottomrule
\end{tabular}
\end{table}

\subsection*{Comparative performance analysis}
To complement the temporal convergence study, we analyse the computational efficiency of the three IMEX schemes by considering the relationship between numerical error, CPU time and the number of time steps for the quantities $v$, $u_e$, and $\FP[5]$. 

\begin{figure}[htbp]
    \centering
    \includegraphics[width=\linewidth,height=0.8\textheight,keepaspectratio]{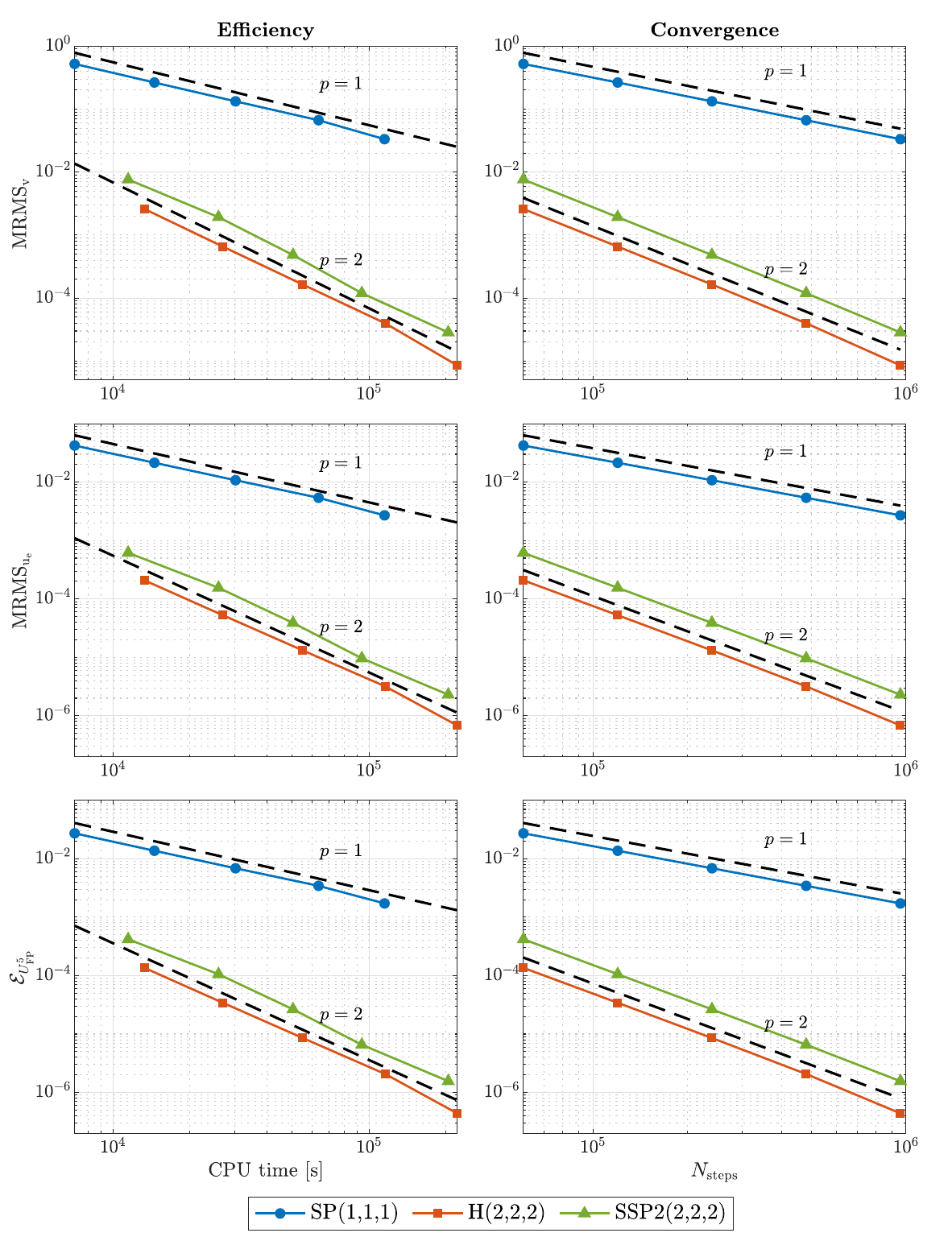}
    \caption{Performance and convergence analysis of the IMEX schemes  for the 2D MEA model. The left column displays the computational efficiency, evaluated as MRMS error versus CPU time, while the right column illustrates the temporal convergence in terms of MRMS error versus the number of time steps ($N_{\mathrm{steps}}$). Results are shown for $v$ (top row),  $u_e$ (middle row), and $\FP[5]$ (bottom row). Solid lines represent the numerical results obtained with SP(1,1,1), H(2,2,2), and SSP2(2,2,2). Dashed lines represent the theoretical slopes for first-order ($p=1$) and second-order ($p=2$) accuracy.}
    \label{Comp_performance}
\end{figure}

The computational times associated with each method for the different time step sizes are summarized in Table \ref{tab:cpu_time}. 
Figure \ref{Comp_performance} presents these results on log-log scales to highlight the asymptotic convergence rates and facilitate a direct comparison of the schemes' efficiency. The error-versus-step curves (right column) confirm the theoretical orders of accuracy ($p=1$ and $p=2$), while the efficiency plots (error versus CPU time, left column) provide a practical measure of performance. Notably, although second-order schemes require more computational effort per single time step, they achieve a target accuracy significantly faster than the first-order SP(1,1,1) scheme by allowing the use of much larger $\Delta t$.

\begin{table}[t]
\setlength{\tabcolsep}{3pt}
\caption{Computational time (hh:min) and accuracy (using $\MRMS_v$) of the IMEX schemes for decreasing values of $\Delta t$ for the 2D MEA model.}
\label{tab:cpu_time}
\begin{center}
\begin{tabular}{l cc cc cc}
\toprule
\multirow{2}{*}{$\Delta t$ [ms]} & \multicolumn{2}{c}{SP(1,1,1)} & \multicolumn{2}{c}{H(2,2,2)} & \multicolumn{2}{c}{SSP2(2,2,2)} \\
\cmidrule(lr){2-3} \cmidrule(lr){4-5} \cmidrule(lr){6-7}
 & \small{CPU time}  & \small{\MRMS$_{v}$} & \small{CPU time}  & \small{\MRMS$_{v}$} & \small{CPU time}  & \small{\MRMS$_{v}$} \\
\midrule
1e-2    & 01:58 & 5.22e-1 & 03:41 & 2.59e-3 & 03:11 & 7.63e-3 \\
5e-3  & 04:02 & 2.65e-1 & 07:27 & 6.57e-4 & 07:10 & 1.93e-3 \\
2.5e-3 & 08:22 & 1.33e-1 & 15:18 & 1.64e-4 & 14:01 & 4.83e-4 \\
1.25e-3 & 17:38 & 6.66e-2 & 32:08 & 3.98e-5 & 26:02 & 1.20e-4 \\
6.25e-4 & 31:55 & 3.33e-2 & 61:20 & 8.55e-6 & 56:36 & 2.85e-5 \\
\bottomrule
\end{tabular}
\end{center}
\end{table}

\subsection*{OS CPU time analysis}
\label{sec:stage_wise_cpu}

To gain deeper insights into the computational efficiency of the investigated methods within the proposed OS framework, we perform a detailed CPU time profiling for each splitting stage (Steps 1 to 4). This stage-based analysis allows us to identify the most computationally demanding parts and map their behaviour to the underlying hiPSC-CM physiological dynamics.

\begin{figure}[t]
    \centering
    \includegraphics[width=0.7\linewidth]{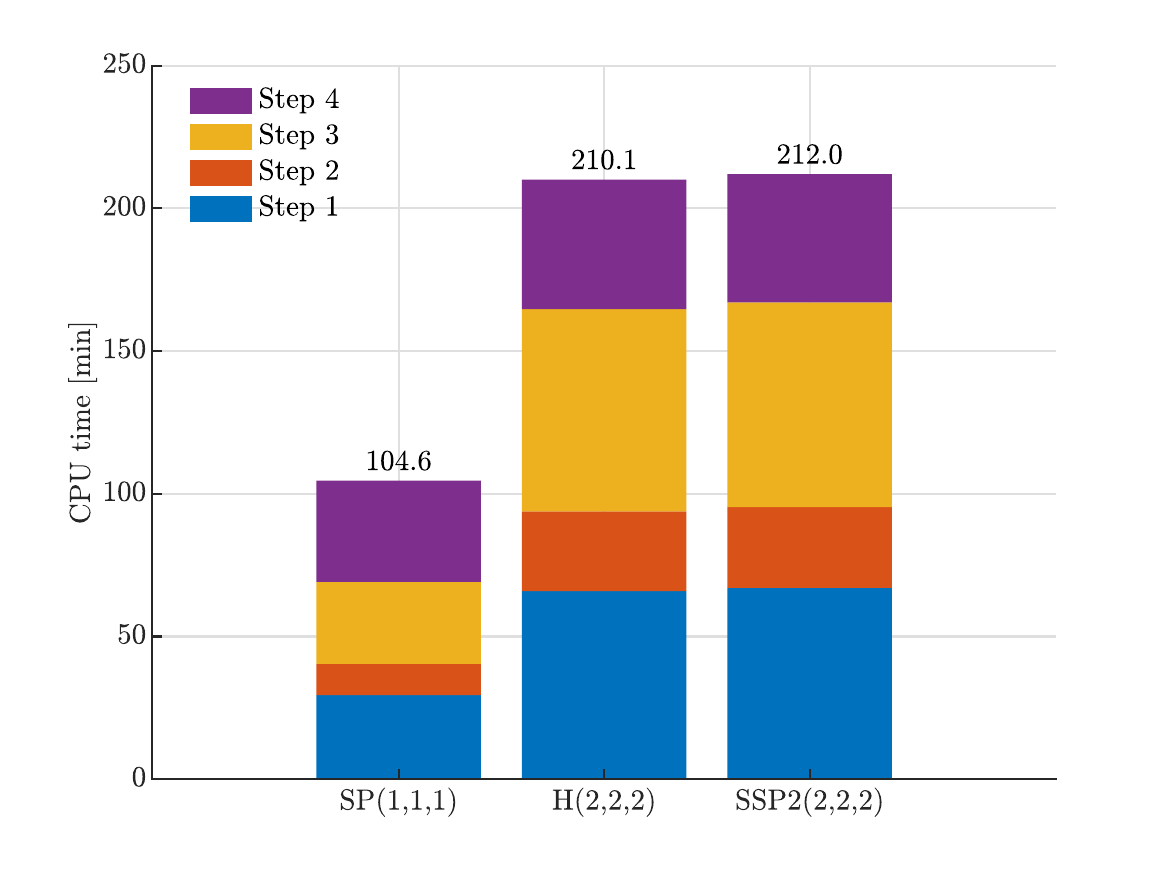}
    \caption{Cumulative computational cost (minutes) for all the IMEX schemes and all the OS stages of the 2D MEA model. CPU times are computed in the case of $\Delta t = \text{1e-2}$ ms.}
    \label{fig:bar_plot}
\end{figure}

Table~\ref{tab:cpu_summary} summarizes both the average per time step CPU times and the cumulative computational costs associated with each OS stage for the considered IMEX schemes. The corresponding cumulative stage-wise computational costs are additionally visualized in Fig.~\ref{fig:bar_plot}.

The results show that the adoption of second-order IMEX schemes approximately doubles the total computational time with respect to the first-order $\mathrm{SP(1,1,1)}$ method. In particular, $\mathrm{H(2,2,2)}$ and $\mathrm{SSP2(2,2,2)}$ schemes exhibit nearly identical computational costs.

\begin{figure}[t]
    \centering
    \includegraphics[width=\linewidth]{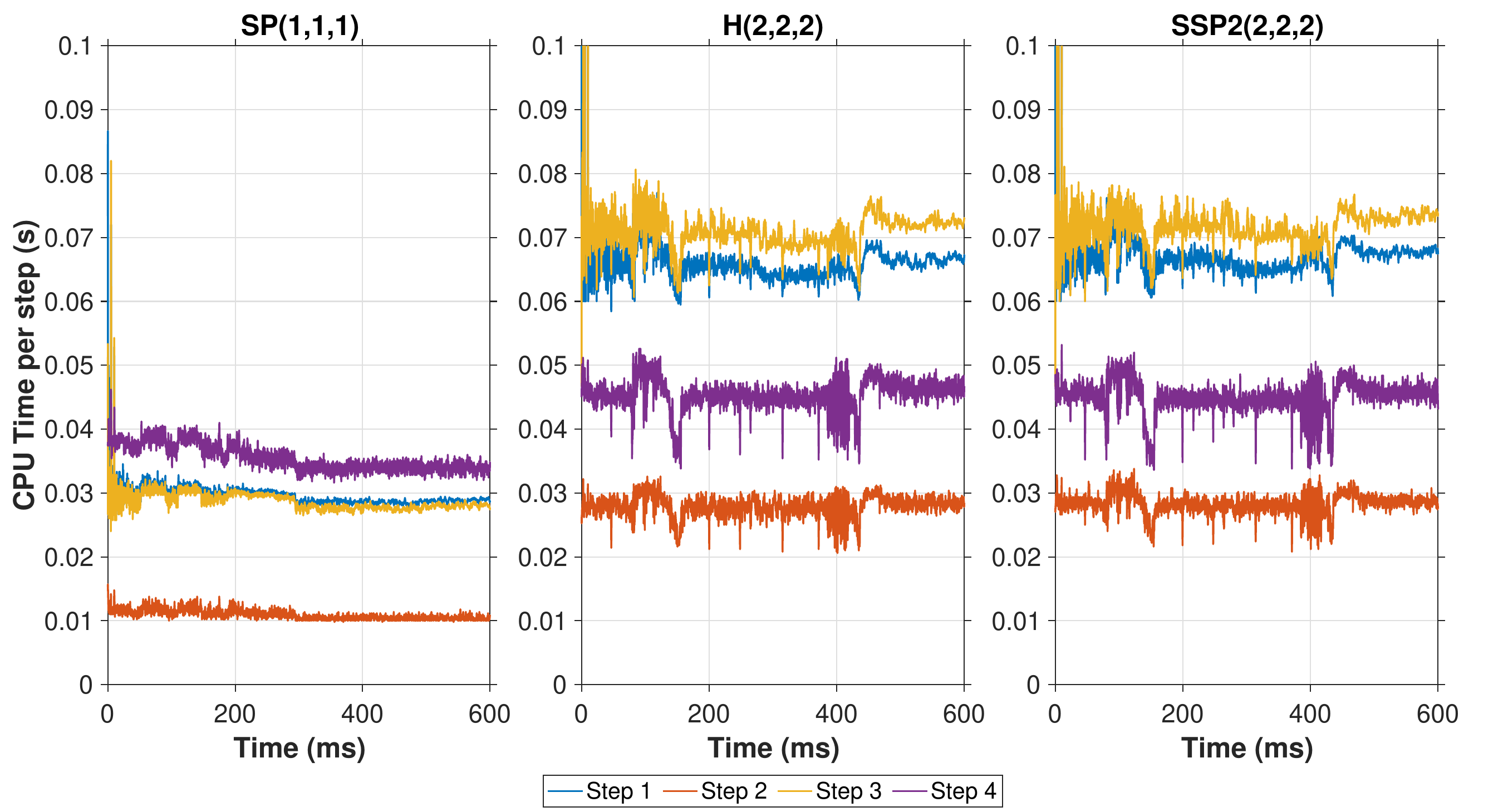}
    \caption{CPU time per time step of the 2D MEA model across the four splitting stages for SP(1,1,1) (left), H(2,2,2) (centre), and SSP2 (2,2,2) (right) over the $600\,\mathrm{ms}$ simulation. Computations are performed using $\Delta t =\text{1e-2}$ ms.}
    \label{cpu_for_time_steps}
\end{figure}

For all considered methods, the dominant contribution to the overall CPU time arises from Steps 1 and 3, corresponding to the solution of the ionic model. Conversely, Steps 2 and 4 remain computationally cheaper, although their cost also increases when moving from first- to second-order schemes. To further investigate the temporal evolution of the computational effort, the CPU time associated with each time step is analysed in Fig.~\ref{cpu_for_time_steps}.

\begin{table}[tb]
\centering
\caption{Average per-timestep execution times and cumulative computational costs associated with the four OS stages for the considered IMEX schemes. CPU times are computed using $\Delta t = \text{1e-2}$ ms.}
\label{tab:cpu_summary}
\setlength{\tabcolsep}{3pt}
\vspace{0.2cm}
\begin{tabular}{lccccccccc}
\toprule
\multirow{3}{*}{IMEX Method}
& \multicolumn{4}{c}{Average CPU Time }
& \multicolumn{4}{c}{Cumulative CPU Time}
& \multirow{3}{*}{Total CPU} \\

& \multicolumn{4}{c}{per time step [s]}
& \multicolumn{4}{c}{per OS Stage [min]}
& \multirow{3}{*}{Time [min]}\\
\cmidrule(lr){2-5}
\cmidrule(lr){6-9}

& Step 1 & Step 2 & Step 3 & Step 4
& Step 1 & Step 2 & Step 3 & Step 4 & \\

\midrule

SP(1,1,1)
& 0.03 & 0.01 & 0.03 & 0.04
& 29.5 & 10.9 & 28.7 & 35.5 & 104.6 \\

H(2,2,2)
& 0.07 & 0.03 & 0.07 & 0.05
& 65.9 & 28.0 & 70.9 & 45.3 & 210.1 \\

SSP2(2,2,2)
& 0.07 & 0.03 & 0.07 & 0.04
& 67.0 & 28.3 & 71.7 & 45.0 & 212.0 \\

\bottomrule
\end{tabular}
\end{table}

A common feature shared by all IMEX schemes is the presence of a sharp localized spike during the initial time steps of the simulation. This transient behaviour correlates with the rapid depolarization phase of the cellular AP, where stiff nonlinear ionic dynamics induce a substantially larger workload. Such an effect is particularly evident in Steps 1 and 3.

After this initial transient, the computational cost progressively stabilizes during the plateau and repolarization phases of the action potential. Within this regime, the first-order $\mathrm{SP(1,1,1)}$ scheme exhibits the lowest CPU time per time step, whereas the two second-order methods maintain almost identical temporal profiles throughout the simulation.

\section{Discussion and conclusions}
\label{sec:conclusion}
In this work, we developed a mathematical and computational framework for the simulation of electrophysiological activity in hiPSC-CM monolayers recorded through MEA devices. The model, based on the bidomain equations and coupled with a detailed electrode representation, enables a consistent reconstruction of FP signals. The numerical results provide insight into both the modelling capabilities of the framework and the performance of the adopted numerical schemes.

From a modelling viewpoint, the simulated signals successfully reproduce the main qualitative features of the electrophysiological response. In particular, the FP traces display the expected morphology. This indicates that the proposed coupling between the Bidomain model and the electrode representation provides a consistent description of the MEA dynamics.

From a computational perspective, the combination of OS and IMEX Runge--Kutta methods proved to be effective for handling the strong nonlinearity and stiffness of the coupled problem. The convergence analysis, performed on both the 0D ionic model and the full 2D MEA configuration, confirmed the theoretical orders of accuracy of all the considered schemes. The observed convergence behaviour was consistent across all the investigated quantities ($v$, $u_e$, and $U_{\mathrm{FP}}$), as well as across the different error metrics adopted in the study.

A key outcome of this work is the clear advantage provided by second-order IMEX schemes. While the first-order scheme is computationally less expensive per time step, its lower accuracy implies that significantly smaller time steps are required to reach a given error tolerance. In contrast, the second-order schemes achieve the same level of accuracy with coarser temporal resolutions, effectively reducing the overall computational effort. This behaviour is clearly reflected in the error versus CPU time plots (Fig.~\ref{fig:Comp_performance_0d} --~\ref{Comp_performance}), where higher-order methods exhibit a better efficiency profile. This behaviour was observed in both the 0D and 2D settings and remained robust across all variables and error measures considered. Quantitatively, in most of the test cases, the higher-order schemes achieved accuracy gains up to 2--3 orders of magnitude with respect to the first-order approach for a fixed time step size.

Overall, the results indicate that second-order IMEX--RK schemes provide an excellent balance between accuracy and computational efficiency. By achieving substantially lower errors for a comparable computational cost, they represent a highly effective choice for accurate large-scale electrophysiological simulations.

\subsection*{Limitations and Future Work} 
Despite the promising results obtained in this study, several aspects deserve further investigation. First, the numerical analysis was restricted to first- and second-order time discretization due to the adoption of a Strang-based operator-splitting framework. Extending the present approach to third- or higher-order schemes, which have shown promising results in cardiac electrophysiology (\cite{cervi2018high,cervi2019fourth,spiteri2016,Spiteri2025improving}), may provide additional gains in accuracy and computational efficiency.

Furthermore, although the numerical experiments confirmed the expected orders of convergence, a rigorous theoretical analysis of the fully coupled nonlinear MEA model is still missing. In particular, convergence, stability, and error estimates for the proposed IMEX operator-splitting discretization remain an open mathematical question. Addressing these issues would provide a stronger theoretical foundation for the methodology and represents an important direction for future research.

\section*{Acknowledgments}
This work was supported by Istituto Nazionale di Alta Matematica (INdAM-GNCS), by the MICROCARD-2 project (ID 101172576) and by the Swiss National Science Foundation, SNSF grants 217025.

\bibliography{refs}  
\end{document}